\documentclass[11pt,english]{article}
\usepackage{mathptmx}
\usepackage[T1]{fontenc}
\usepackage[latin9]{inputenc}
\usepackage[a4paper]{geometry}
\usepackage{float}
\usepackage{amsmath}
\usepackage{amssymb}
\usepackage{graphicx}

\makeatletter

\providecommand{\tabularnewline}{\\}

\makeatother

\usepackage{babel}
\begin{document}
\title{A study of spectral element method for elliptic interface problems
with nonsmooth solutions in $\mathbb{R}^{2}$}
\author{N. Kishore Kumar\thanks{BITS-Pilani Hyderabad Campus, Hyderabad , India; Email: naraparaju@hyderabad.bits-pilani.ac.in }\thanks{This research work is supported by National Board of Higher Mathematics,
DAE, India. } , Pankaj Biswas\thanks{National Institute of Technology Silchar, Silchar, India, Email:pankaj@math.nits.ac.in }
~and B. Seshadri Reddy\thanks{BITS-Pilani Hyderabad Campus, Hyderabad}}
\date{~}
\maketitle
\begin{abstract}
The solution of the elliptic partial differential equation has interface
singularity at the points which are either the intersections of interfaces
or the intersections of interfaces with the boundary of the domain.
The singularities that arises in the elliptic interface problems are
very complex. In this article we propose an exponentially accurate
nonconforming spectral element method for these problems based on
\cite{PRAVIKISH,KISHORENAGA2}. A geometric mesh is used in the neighbourhood
of the singularities and the auxiliary map of the form $z=\ln\xi$
is introduced to remove the singularities. The method is essentially
a least-squares method and the solution can be obtained by solving
the normal equations using the preconditioned conjugate gradient method
(PCGM) without computing the mass and stiffness matrices. Numerical
examples are presented to show the exponential accuracy of the method. 
\end{abstract}
\textbf{}

\section{Introduction}

An interface problem is a special case of an elliptic differential
equation with discontinuous coefficients. Such interface problems
arise in different situations, for example, in heat conduction or
in elasticity problems whose solution domains are composed of several
different materials. There are different kinds of elliptic interface
problems: the interface problems with smooth interfaces, the interface
problems with nonsmooth solutions etc. When the interface is smooth
enough the solution of the interface problem is also very smooth in
the individual regions but global regularity is low, i.e the solution
$u\in H^{1}(\Omega)$ and $u\notin H^{k}(\Omega)$ for $k\geq2.$
This case has been widely addressed in the literature using finite
element methods \cite{Babuska3,bareli,brambthkin}, immersed interface
methods \cite{liito} and least-squares methods \cite{caogunz} etc.
For further information on this problem and existing numerical approaches
in the literature, refer to \cite{KISHORENAGA2}. In this article
we consider the interface problems with nonsmooth solutions. 

In the solution of the elliptic boundary value problems, singularities
may occur when the boundary is not smooth \cite{grisward} or when
the boundary is smooth yet one or more data are not smooth. The second
type of singularity typically arises in interface problems. The singularities
that arises in the interface problems are very complex. The solution
of the elliptic differential equation has interface singularity at
the points which are either the intersections of interfaces or the
intersections of interfaces with the boundary of the domain. The solution
also has singular behavior at the points where the interfaces crosses
each other. The interface singularity at the crossing of interfaces
is very strong. 

The singularities in interface problems has been studied by Kellogg
(considered the interface problem for Poisson equation) in \cite{kellog 1}.
In \cite{kellogg2} Kellogg had studied the Poisson equation with
intersecting interfaces. The complexity depends on the structure of
the eigenvalues of Sturm-Liouville problems corresponding to the singularities.
The elliptic interface problems with singularities also has been studied
in \cite{kellog3,Nicaise,Petzoldt}.

The conventional numerical approaches (the finite difference as well
as finite element) may fail to provide any practical engineering accuracy
at a reasonable cost. In \cite{Babuska3} Babuska studied the interface
problem in the frame work of finite element method. The rates of convergence
are algebraic for the $h-$version and $p-$version of the finite
element method. The mesh refinements techniques gives reasonably good
results but they require longer computing time and also cannot give
acceptable result when the singularities are very strong. The method
of auxiliary map has been introduced in \cite{ohbabu} for the interface
problems by Oh and Babuska in the framework of $p-$ version of FEM
(originally introduced in \cite{lucas} for elliptic problems containing
singularities as MAM). With a proper choice of auxiliary mappings
this method gives better results than the mesh refinements when the
interface singularities are very strong. An optimal choice of the
auxiliary mappings requires a prior knowledge of the structure of
the interface singularities at the singular points. 

In \cite{babguo} an exponentially accurate method ($hp$ finite element)
has been proposed by Babuska and Guo for the elliptic problems with
analytic data on the nonsmooth domains like the domains with cracks,
re-entrant corners. Geometric mesh has been considered near the corners
to resolve the singularities in the solution. They have studied the
regularity of the solution in the framework of weighted Sobolev space
$H_{\beta}^{k,2}(\Omega)$ and the countably normed space $B_{\beta}^{l}(\Omega)$.
In \cite{babguo2} Babuska and Guo have analyzed the regularity of
the interface problem in terms of countably normed spaces. In \cite{guooh}
Guo and H. S. Oh analyzed the $hp$ version of the finite element
method for problems with interfaces. They have used geometric mesh
near the singularities and shown the exponential accuracy. Geometric
mesh together with the auxiliary mapping technique gives better results
even if the singularities are extremely severe. They have also presented
the theoretical results for interface problems. The theoretical results
and numerical scheme of this version can also be applied to general
elliptic equations and systems, including elasticity problems with
homogeneous and non-homogeneous materials. 

In \cite{honhuang} H. Hon and Z. Huang introduced the direct method
of lines for numerical solution of interface problems. The interface
problem is reduced to variational-differential (V-D) problem on semi-infinite
strip in $\rho$ and $\phi$ variables by using a suitable transformation
of coordinates. Furthermore, the V-D problem is discretized respect
with the variable $\phi$ and solved numerically. This method requires
no prior knowledge of the constructure of the singularity at the singular
point. 

In \cite{PRAVIKISH,KISHOREPRAVIR,KISHORENAGA} P. Dutt et. al. proposed
an exponentially accurate nonconforming $hp$/spectral element method
to solve general elliptic boundary value problems with mixed Neumann
and Dirichlet boundary conditions on non-smooth domains. In \cite{KISHORENAGA2},
a spectral element method for elliptic interface problems with smooth
interfaces has been introduced and this has been extended to the elasticity
interface problems in \cite{kishore}. Blending elements have been
used to completely resolve the interface and higher order approximation
has been used. 

In this article we propose a nonconforming spectral element method
for elliptic interface problems with singularities based on \cite{PRAVIKISH,KISHOREPRAVIR,KISHORENAGA}.
A geometric mesh is used in the neighbourhood of the vertices and
the auxiliary map of the form $z=\ln\xi$ is introduced to remove
the singularities at the corners, which was first introduced by Kondratiev
in \cite{KONDR}. In the remaining part of the domain usual Cartesian
coordinate system is used. The proposed method is essentially a least-squares
method. 

In the least-squares formulation of the method, a solution is sought
which minimizes the sum of the squares of a weighted squared norms
of the residuals in the partial differential equation and the sum
of the squares of the residuals in the boundary conditions in fractional
Sobolev norms and the sum of the squares of the jumps in the value
and its normal derivatives of the function across the interface in
appropriate fractional Sobolev norms and enforce the continuity along
the inter element boundaries by adding a term which measures the sum
of the squares of the jump in the function and its derivatives in
fractional Sobolev norms.

The spectral element functions are nonconforming. The solution can
be obtained by solving the normal equations using the preconditioned
conjugate gradient method (PCGM) without computing the mass and stiffness
matrices \cite{KISHOREPRAVIR,tomar}. An efficient preconditioner
is used for the method which is a block diagonal matrix, where each
diagonal block corresponds to an element \cite{pravirpankaj}. The
condition number of the preconditioner is $O\left(\ln W\right)^{2},$
where $W$ is the degree of the approximating polynomial. Let $N$
denote the number of layers in the geometric mesh such that $W$ is
proportional to $N$. Then the method requires $O(W\,\ln W)$ iterations
of the PCGM to obtain the solution to exponential accuracy.

Here we define some Sobolev norms which are used in this article.
Denote by $H^{m}(\Omega)$ the Sobolev space of functions with square
integrable derivatives of integer order $\leq m$ on $\Omega$ (a
domain) furnished with the norm
\[
\parallel u\parallel_{H^{m}(\Omega)}^{2}=\sum_{\mid\alpha\mid\leq m}\parallel D^{\alpha}u\parallel_{L^{2}(\Omega)}^{2}.
\]

Further, let 
\begin{eqnarray*}
 &  & \left\Vert u\right\Vert _{s,I}^{2}=\int_{I}u^{2}(x)dx+\int_{I}\int_{I}\frac{\left|u(x)-u(x^{\prime})\right|^{2}}{\left|x-x^{\prime}\right|^{1+2s}}dxdx^{\prime}
\end{eqnarray*}
 denote the fractional Sobolev norm of order $s,$ where $0<s<1.$
Here $I$ denotes an interval contained in $\mathbb{R}.$ 

For the definitions of the other function spaces which appears in
this article, refer to \cite{babguo,grisward,guooh}. Throughout the
article we use $x=(x_{1},x_{2})$ to represent a point on $\mathbb{R}^{2}$
(in Cartesian coordinate system).

The contents of this paper are organized as follows: In Section 2
the elliptic interface problem is defined. Discretization of the domain
and local transformation are described in Section 3. In Section 4,
the numerical scheme has been derived. In Section 5 the computational
results are provided for few test problems. 

\section{Interface Problem}

In this section we state define the elliptic interface problem. For
the convenience of the reader, we consider the polygonal domain as
shown in the figure 1 for defining the interface problem, discretization
and deriving the numerical scheme. The numerical method is also applicable
for general polygonal domains with more number of vertices. 

\begin{figure}[H]
~~~~~~~~~~~~~~~~~~~~~~~~~~~~~~~~~~~~~~~~~~~~\includegraphics[scale=0.8]{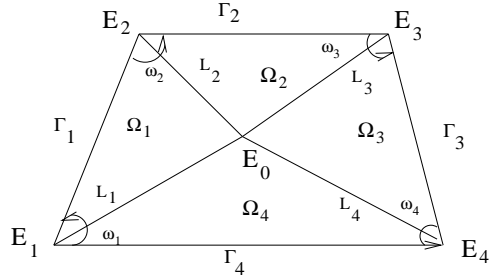}

\caption{Polygonal domain with interfaces}
\end{figure}

Consider the polygonal domain $\Omega$ in $\mathbb{R}^{2}$ with
boundary $\partial\Omega=\Gamma$ as shown in the Fig. 1. Let $E_{i},i=1,2,3,4$
be the vertices of the domain. Let $\varGamma={\displaystyle \cup_{i=1}^{4}}\Gamma_{i},$
where $\Gamma_{i}$ be the open edge of $\partial\Omega$ connecting
$E_{i}$ and $E_{i+1}$. By $\omega_{i}$ we denote the measure of
the interior angle of $\Omega$ at $E_{i}.$ Without loss of generality,
we will assume all interfaces meet at only one point $E_{0}$ as shown
in Fig. 1. Let $(r,\theta)=(r_{0},\theta_{0})$ be the polar co-ordinates
at the point $E_{0}$ and suppose $\Omega$ is decomposed into four
subdomains $\Omega_{1},\Omega_{2},\Omega_{3},\Omega_{4}$ so that
$\overline{\Omega}_{k}\cap\overline{\Omega}_{k-1}$ is a straight
line interface $L_{k},$ $L_{k}=\{(r,\theta):\theta=\Theta_{k},0\leq r\leq R_{k}\},$
for $k=1,2,3,4.$ \medskip{}

\subsubsection*{Elliptic Interface problem}

Let us consider the following interface problem:
\begin{eqnarray}
 &  & \mathcal{L}u=-\nabla.(p\nabla u)=f\,\,\textrm{in}\,\cup\Omega_{i}\nonumber \\
 &  & u=0\,\,\textrm{on}\,\,\Gamma_{D}=\cup_{i\in\mathcal{D}}\overline{\Gamma}_{i}\nonumber \\
 &  & \frac{\partial u}{\partial n}=g=G^{N}\mid_{\Gamma_{N}}\,\,\textrm{on}\,\,\Gamma_{N}=\cup_{i\in\mathcal{N}}\overline{\Gamma}_{i}\label{eq:1}
\end{eqnarray}
 where $\Gamma_{D}\cup\Gamma_{N}=\partial\Omega,\mathcal{D}\cup\mathcal{N}=\{1,2,..,4\},\mathcal{D}\cap\mathcal{N}=\emptyset,n=(n_{1},n_{2})$
is the unit normal vector on $\Gamma_{N}$ and the coefficients are
piecewise constants:
\begin{eqnarray}
 &  & p=\begin{cases}
p_{1} & \textrm{in}\,\,\Omega_{1}\\
p_{2} & \textrm{in}\,\,\Omega_{2}\\
. & .\\
p_{4} & \textrm{in}\,\,\Omega_{4}
\end{cases}.\label{eq:2}
\end{eqnarray}

Assume that the interface conditions are satisfied. That is, on $L_{k},$
$1\leq k\leq4,$ $u$ satisfies
\begin{eqnarray}
 &  & u(r,\Theta_{k}-0)=u(r,\Theta_{k}+0)\nonumber \\
 &  & p_{k-1}\frac{\partial u}{\partial n}(r,\Theta_{k}-0)=p_{k}\frac{\partial u}{\partial n}(r,\Theta_{k}+0)\label{eq:3}
\end{eqnarray}
 where $n=(n_{1},n_{2})$ is a unit normal vector to the interfaces
$L_{k}.$ 

The asymptotic expansion, uniqueness and regularity of the solution
of the above problem (1) - (3) has been discussed in detail in \cite{guooh}.
It has been shown that the solution has $r^{\lambda}$ type of singularity
near the points $E_{i}$ which is similar to the singularity in the
solution of elliptic problems on nonsmooth domains like domains with
cracks and reentrant corners. But the strength of the singularity
is strong in the elliptic interface problems.

\section{Discretization and Stability}

\begin{figure}[H]

~~~~~~\includegraphics[scale=0.75]{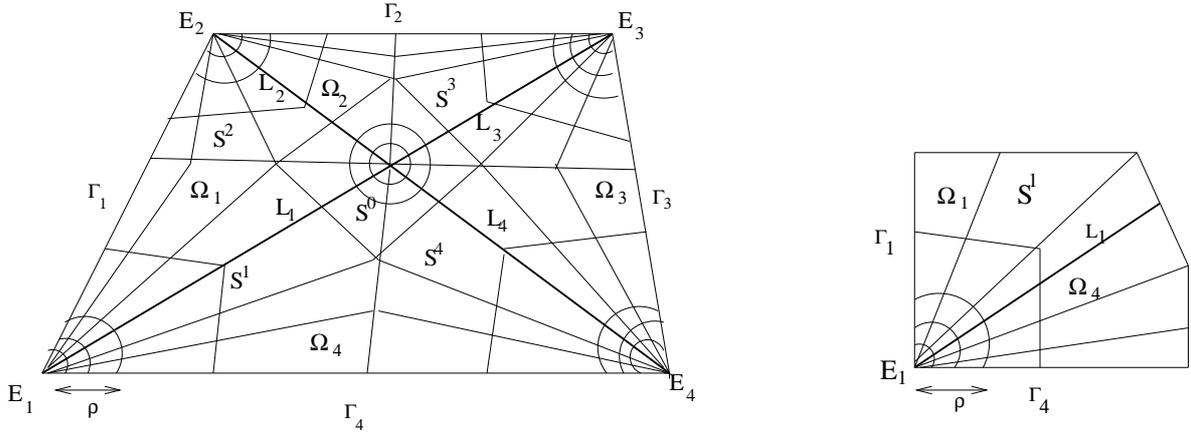}

\caption{Discretization of the domain}

\end{figure}

Discretize the polygonal domain $\Omega$ into $5$ non-overlapping
polygonal subdomains $S^{0},S^{1},\ldots,S^{4}.$ Here each $S^{k}\,\,\textrm{for}\,\,k=1,2,3,4$
contains the vertex $E_{k}$ only and contains a part of the interface
$L_{k}$ and $S^{0}$ contains the point $E_{0}$ and contains a part
of all the interfaces $L_{k}$ as shown in Fig. $\negmedspace$2.
Each subdomain $S^{k},\,k\neq0$ is a subset of union of two subdomains
$\Omega_{i}$ and $\Omega_{j}$ for some $i,j.$ For example $S^{1}\subset\Omega_{1}\cup\Omega_{4}$
as shown in Fig. 2. 

Let $S^{k}=\left\{ \Omega_{i,j}^{k}:\,j=1,2,..,J_{k},\,i=1,2,\ldots,I_{k}\right\} $
be a partition of $S^{k},k=1,2,3,4$ where $J_{k}$ and $I_{k}$ are
integers. Let $I_{k}$ be bounded and constant for all $k=1,2,3,4.$
Let $S^{0}=\left\{ \Omega_{i,j}^{0}:\,j=1,2,..,J_{0},\,i=1,\right.$
$\left.2,\ldots,I_{0}=8\right\} $ be a partition of $S^{0}$ such
that subdomain division matches on the interface. Let $\left(r_{k},\theta_{k}\right)$
denote polar coordinates with center at $E_{k}.$ 

Since the solution of the interface problem has singular behavior
at $E_{i},i=1,2,3,4$ where the interface intersects the boundary
and also at $E_{0}$ where the interfaces meet each other, we consider
the geometric mesh and use the auxiliary mapping near each point $E_{i}.$
The description of the geometric mesh and the auxiliary mapping is
given below. 

\subsubsection*{Geometric mesh near $E_{k},k\protect\neq0$}

Let $\left\{ \psi_{i}^{k}\right\} _{i=1,\ldots,I_{k}+1}$ be an increasing
sequence of points such that $\psi_{1}^{k}=\psi_{l}^{k}$ and $\psi_{I_{k}+1}^{k}=\psi_{u}^{k}.$
Let $\psi_{i}^{k}$ meet with interface for some $i=I$. That is,
$\psi_{I}^{k}$ matches with the interface $L_{k}$ and hence separates
$\Omega_{i}$ and $\Omega_{j}$ in $S^{k}.$ Let $\Delta\psi_{i}^{k}=\psi_{i+1}^{k}-\psi_{i}^{k}.$
Choose these points so that 
\begin{eqnarray*}
\max_{k}\left(\max_{i}\Delta\psi_{i}^{k}\right) & \leq & \lambda\min_{k}\left(\min_{i}\Delta\psi_{i}^{k}\right)
\end{eqnarray*}
 for some constant $\lambda.$ 

Let $\Pi^{k}=\{(x_{1},x_{2}):0<r_{k}<\rho\}\subseteq S^{k}$ be a
sector with sides $\varGamma_{k}$ and $\varGamma_{k+1}.$ Now choose
a geometric mesh with $N$ layers in $\Pi^{k}$ with a geometric ratio
$q_{k}\left(0<q_{k}<1\right).$ Let $\sigma_{j}^{k}=\rho\left(q_{k}\right)^{N+1-j}\;\textrm{for}\:2\leq j\leq N+1$
and $\sigma_{1}^{k}=0.$ 

Let
\begin{eqnarray*}
\Omega_{i,j}^{k} & = & \left\{ \left(x_{1},x_{2}\right):\sigma_{j}^{k}<r_{k}<\sigma_{j+1}^{k},\psi_{i}^{k}<\theta_{k}<\psi_{i+1}^{k}\right\} ,
\end{eqnarray*}
$\mathrm{\textrm{for}}\:1\leq i\leq I_{k},1\leq j\leq N.$

Since $S^{k}$ contains a part of the interface $L_{k}$ and $\psi_{I}^{k}$
meet with it, the elements $\Omega_{I,j}^{k}$ and $\Omega_{I-1,j}^{k}$
have the common edge which lies on the interface. For example, the
elements $\Omega_{I,j}^{1}\subset\Omega_{1}$ and $\Omega_{I-1,j}^{1}\subset\Omega_{4}$
in $S^{1}$ have the common edge on the interface $L_{1}.$ 

\subsubsection*{Geometric mesh near $E_{0}$}

Let $\left\{ \psi_{i}^{0}\right\} _{i=1,\ldots,9}$ be an increasing
sequence of points such that $\psi_{1}^{0}=0$ and $\psi_{9}^{0}=2\pi$
and $\psi_{i}^{0}$ for some $i$ meet with interfaces $L_{k}$. 

Let $\Pi^{0}=\{(x_{1},x_{2}):0<r_{0}<\rho\}\subseteq S^{0}$ be a
circular region around $E_{0}$. Now choose a geometric mesh with
$N$ layers in $\Pi^{0}$ with a geometric ratio $q_{0}\left(0<q_{0}<1\right).$
Let $\sigma_{j}^{0}=\rho\left(q_{0}\right)^{N+1-j}\;\textrm{for}\:2\leq j\leq N+1$
and $\sigma_{1}^{0}=0.$ 

Let
\begin{eqnarray*}
\Omega_{i,j}^{0} & = & \left\{ \left(x_{1},x_{2}\right):\sigma_{j}^{0}<r_{0}<\sigma_{j+1}^{0},\psi_{i}^{0}<\theta_{0}<\psi_{i+1}^{0}\right\} ,
\end{eqnarray*}
$\mathrm{\textrm{for}}\:1\leq i\leq8,1\leq j\leq N.$

\subsubsection*{In the remaining part of $S^{k}$}

In the remaining part of $S^{k},1\leq k\leq4,$ we retain the Cartesian
coordinate system $(x_{1},x_{2})$ i.e., in $\Omega_{i,j}^{k}$ for
$1\leq i\leq I_{k},N<j\leq J_{k}.$

Let
\begin{eqnarray*}
 &  & \Omega^{1}=\left\{ \Omega_{i,j}^{k}:1\leq i\leq I_{k},N<j\leq J_{k},1\leq k\leq4\right\} .
\end{eqnarray*}

Similarly we retain the Cartesian coordinate system $(x_{1},x_{2})$
in the remaining part of $S^{0}.$

Let 
\begin{eqnarray*}
 &  & \Omega^{0}=\left\{ \Omega_{i,j}^{0}:1\leq i\leq8,N<j\leq J_{k}\right\} .
\end{eqnarray*}

Here for $i=1,8$, $\Omega_{i,j}^{0}\subseteq\Omega_{3};$~~$i=2,3,$~~$\Omega_{i,j}^{0}\subseteq\Omega_{2};$~~~$i=4,5,$~~$\Omega_{i,j}^{0}\subseteq\Omega_{1};$~~$i=6,7,$~~$\Omega_{i,j}^{0}\subseteq\Omega_{4}.$
For $j>N,$ $\Omega_{1,j}^{0}$ and $\Omega_{2,j}^{0}$ have a common
edge which lies on $L_{3}.$ Similarly, the elements $\Omega_{1,j}^{0}\,\&\,\Omega_{2,j}^{0}$
$\Omega_{3,j}^{0}\,\&\,\Omega_{4,j}^{0},$ $\Omega_{5,j}^{0}$ \&
$\Omega_{6,j}^{0}$ and $\Omega_{7,j}^{0}$ \& $\Omega_{8,j}^{0}$
have the common edges which lies on $L_{3},L_{2},L_{1}$ and $L_{4}$
respectively. 

\subsubsection*{Auxiliary Mapping }

Now let $\tau_{k}=\ln\,r_{k}$ in $\{(x_{1},x_{2}):0<r_{k}<\rho\}\subseteq S^{k},$
$0\leq k\leq4.$ Define $\zeta_{j}^{k}=\ln\,\sigma_{j}^{k}$ for $1\leq j\leq N+1$.
Here $\zeta_{1}^{k}=-\infty.$ Define
\begin{eqnarray*}
 &  & \widetilde{\Omega}_{i,j}^{k}=\left\{ \left(\tau_{k},\theta_{k}\right):\;\zeta_{j}^{k}<\tau_{k}<\zeta_{j+1}^{k},\;\psi_{i}^{k}<\theta_{k}<\psi_{i+1}^{k}\right\} ,
\end{eqnarray*}
 for $1\leq i\leq I_{k},1\leq j\leq N.$ Hence the geometric mesh
$\Omega_{i,j}^{k},\,2\leq j\leq N$ becomes a quasi-uniform mesh in
modified polar coordinates (Fig. 3). However, $\widetilde{\Omega}_{i,1}^{k}$
is a semi-infinite strip.

\begin{figure}[H]
~~~~~~~~~~~~~~~~\includegraphics[scale=0.75]{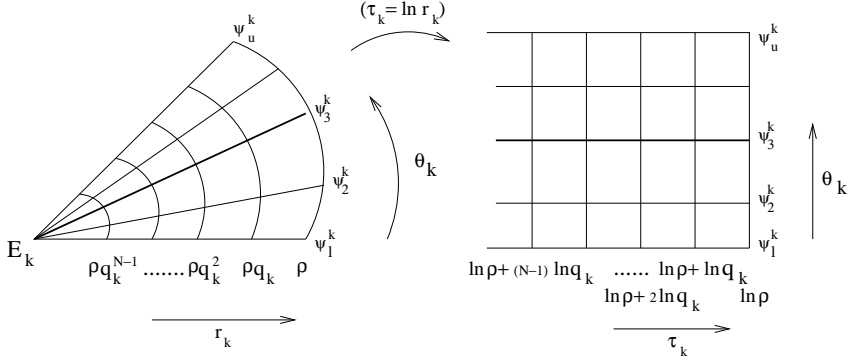}

\caption{Quasi uniform mesh in $\tau_{k}$ and $\theta_{k}$ coordinates near
$E_{k},k\protect\neq0.$}
\end{figure}

\subsubsection*{Approximation }

\hspace*{0.6cm}The nonconforming spectral element functions are sum
of tensor products of polynomials of degree $W_{j},\;1\le W_{j}\le W$
in their respective modified polar coordinates (\ref{eq:a1}) in $\widetilde{\Omega}_{i,j}^{k}$
for $0\le k\le4,1\leq i\leq I_{k},2\leq j\leq N.$ In the infinite
sector i.e., in $\widetilde{\Omega}_{i,1}^{k},$ the solution is approximated
by a constant which is the value of the function $u$ at the corresponding
point $E_{k}$. The constant value is computed by treating it as a
common boundary value during the numerical computation.

Let $u_{i,1}^{k}(\tau_{k},\theta_{k})=h_{k},$ a constant on $\widetilde{\Omega}_{i,1}^{k}$.
Define the spectral element function 
\begin{equation}
u_{i,j}^{k}(\tau_{k},\theta_{k})=\sum_{r=0}^{W_{j}}\sum_{s=0}^{W_{j}}g_{r,s}\,\tau_{k}^{r}\,\theta_{k}^{s},\label{eq:a1}
\end{equation}
on $\widetilde{\Omega}_{i,j}^{k}$ for $1\leq i\leq I_{k},2\leq j\leq N,0\leq k\leq4.$
Here $1\leq W_{j}\leq W.$

Moreover there is an analytic mapping $M_{i,j}^{k}$ from the master
square $S=(-1,1)^{2}$ to the elements $\Omega_{i,j}^{k}$ in $\Omega^{0}$
and $\Omega^{1}.$ Define
\begin{equation}
u_{i,j}^{k}(M_{i,j}^{k}(\xi,\eta))=\sum_{r=0}^{W}\sum_{s=0}^{W}g_{r,s}\,\xi^{r}\,\eta^{s}.\label{eq:a2}
\end{equation}

\section{Numerical Scheme}

Here we describe the numerical formulation. This numerical method
is essentially a least-squares formulation based on \cite{KISHOREPRAVIR,KISHORENAGA2}. 

As defined in Section 3, $\widetilde{\Omega}_{i,j}^{k}$ is the image
of $\Omega_{i,j}^{k}$ in $(\tau_{k},\theta_{k})$ coordinates. Let
$\mathcal{L}_{i,j}^{k}$ be the operator defined by $\mathcal{L}_{i,j}^{k}u=r_{k}^{2}\,\mathcal{L}u.$
Then the operator $\mathcal{\widetilde{L}}_{i,j}^{k}$ in the transformed
coordinates $\tau_{k}\,\,\textrm{and}\,\,\theta_{k}$ is given by
\begin{eqnarray*}
 &  & \tilde{\mathcal{L}}_{i,j}^{k}u=-p\left(\frac{\partial^{2}u}{\partial\tau_{k}^{2}}+\frac{\partial^{2}u}{\partial\theta_{k}^{2}}\right).
\end{eqnarray*}
 Where $p$ takes different values based on $i$ value as explained
in Section 3. 

Next, let the vertex $E_{k}=\left(x_{1}^{k},x_{2}^{k}\right)$ and
\begin{eqnarray*}
 &  & F_{i,j}^{k}\left(\tau_{k},\theta_{k}\right)=e^{2\tau_{k}}f\left(x_{1}^{k}+e^{\tau_{k}}\cos\theta_{k},x_{2}^{k}+e^{\tau_{k}}\sin\theta_{k}\right)
\end{eqnarray*}
in $\tilde{\Omega}_{i,j}^{k}$ for $0\leq k\leq4,$ $2\leq j\leq N,$
$1\leq i\leq I_{k}.$

\begin{figure}[H]
~~~~~~~~~~~~~~~~~~~~~~~~~~~~~~~~~~~~~~~~~~~~~~~~\includegraphics[scale=0.75]{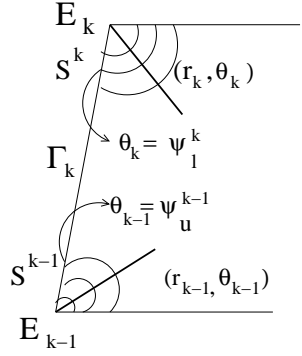}

\caption{Edge $\Gamma_{k}$ common to $\Pi^{k-1}$ and $\Pi^{k}$}
\end{figure}

Consider the boundary $\frac{\partial u}{\partial n}=g$ on $\Gamma_{k}\cap\partial\Pi^{k}$
for $k\in\mathcal{N}$(see Fig. 4). Let
\[
l_{1}^{k}(\tau_{k})=\frac{\partial u^{k}}{\partial n}=e^{\tau_{k}}\,g(x_{1}^{k}+e^{\tau_{k}}cos(\psi_{l}^{k})\,,\,x_{2}^{k}+e^{\tau_{k}}sin(\psi_{l}^{k})).
\]
 Consider $\frac{\partial u}{\partial n}=g$ for $k\in\mathcal{N}$
on $\Gamma_{k}\cap\partial\Pi^{k-1}$(look at Fig. 4). Define 
\[
l_{2}^{k}(\tau_{k-1})=\frac{\partial u^{k}}{\partial n}=e^{\tau_{k-1}}g(x_{1}^{k-1}+e^{\tau_{k-1}}cos(\psi_{u}^{k-1})\,,\,x_{2}^{k-1}+e^{\tau_{k-1}}sin(\psi_{u}^{k-1})).
\]

\begin{figure}[H]
~~~~~~~~~~~~~~~~~~~~~~~~~~~\includegraphics[scale=0.55]{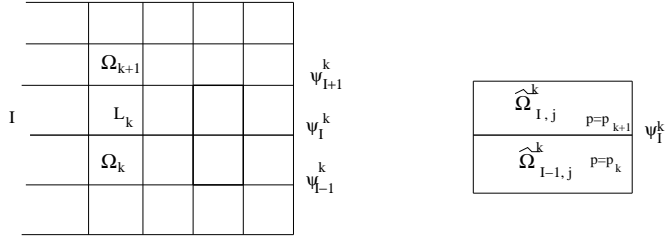}

\caption{Elements with interface as common edge}

\end{figure}

As described in section 3, $\psi_{I}^{k}$ matches with the interface
$L_{k},$ the elements $\Omega_{I,j}^{k}$ and $\Omega_{I-1,j}^{k}$
have the common edge $\gamma_{s}$. Let $\tilde{\gamma}_{s}$ be the
image of $\gamma_{s}$ in $\tau_{k}\,\,\textrm{and}\,\,\theta_{k}$
coordinates and therefore $\tilde{\gamma}_{s}$ is the common edge
of $\tilde{\Omega}_{I,j}^{k}$ and $\tilde{\Omega}_{I-1,j}^{k}$ which
lies on $L_{k}.$ We define the jump in the solution across $\gamma_{s}\subseteq L_{k}$ 

\begin{eqnarray*}
 &  & \left\Vert [u^{k}]\right\Vert _{\frac{3}{2},\tilde{\gamma}_{s}}^{2}=\left\Vert u_{I,j}^{k}(\tau_{k},\psi_{I}^{k})-u_{I-1,j}^{k}(\tau_{k},\psi_{I}^{k})\right\Vert _{_{0,\tilde{\gamma}_{s}}}^{2}+\left\Vert \frac{\partial u_{I,j}^{k}}{\partial\tau_{k}}(\tau_{k},\psi_{I}^{k})-\frac{\partial u_{I-1,j}^{k}}{\partial\tau_{k}}(\tau_{k},\psi_{I}^{k})\right\Vert _{\frac{1}{2},\tilde{\gamma}_{s}}^{2}.
\end{eqnarray*}

Now we define the jump across the normal derivative across the interface
\begin{eqnarray*}
\left\Vert \left[p\frac{\partial u^{k}}{\partial\theta_{k}}\right]\right\Vert _{\frac{1}{2},\tilde{\gamma}_{s}}^{2}=\left\Vert p_{k+1}\frac{\partial u_{I,j}^{k}}{\partial\theta_{k}}(\tau_{k},\psi_{I}^{k})-p_{k}\frac{\partial u_{I-1,j}^{k}}{\partial\theta_{k}}(\tau_{k},\psi_{I}^{k})\right\Vert _{\frac{1}{2},\tilde{\gamma}_{s}}^{2} & .
\end{eqnarray*}

In similar way, we define the term which measures the sum of the squares
of the jump in $u$ and its derivatives with respect to $\tau_{k}$
and $\theta_{k}$ in appropriate Sobolev norms along the inter-element
boundaries. 

In the remaining part of the domain, i.e on $\Omega^{1}$ and $\Omega^{0}$
the solution is smooth. The residue in the equation and jumps across
the interfaces and inter element boundaries and residue at the boundary
has been described in detail in \cite{KISHORENAGA2}. Here we define
the functional near the singularities and in the interior. 

Let $\gamma_{s}\subseteq\bar{\varPi}^{k}$ and $d(E_{k},\gamma_{s})=inf_{_{x\in\gamma_{s}}}\left\{ \textrm{distance}(E_{k},x)\right\} .$
Choose $\alpha_{k}=1-\beta_{k}$ as defined in \cite{PRAVIKISH}.
Let $\mathcal{F}_{u}=\left\{ \left\{ u_{i,j}^{k}(\tau_{k},\theta_{k})\right\} _{i,j,k},\left\{ u_{i,j}^{k}(\xi,\eta)\right\} _{i,j,k}\right\} \in\Pi^{N,W},$
the space of spectral element functions. Define $a_{k}=u(E_{k})$. 

Define the functional
\begin{eqnarray}
\mathfrak{\mathcal{\mathfrak{r}}}_{vertices}^{^{N,W}}(\mathcal{F}_{u}) & = & \sum_{k=0}^{4}\sum_{j=2}^{N}\sum_{i=1}^{I_{k}}(\rho\mu_{k}^{N+1-j})^{-2\alpha_{k}}\left\Vert (\widetilde{\mathcal{L}}_{i,j}^{k})u_{i,j}^{k}(\tau_{k},\theta_{k})-F_{i,j}^{k}\left(\tau_{k},\theta_{k}\right)\right\Vert _{_{0,\tilde{\Omega}_{i,j}^{k}}}^{2}\nonumber \\
 & + & \sum_{k=0}^{4}\sum_{{\gamma_{s}\subseteq\varPi^{k},\gamma_{s}\varsubsetneq L_{k}\atop \mu(\tilde{\gamma_{s}})<\infty}}\!\!\!d(E_{k},\gamma_{s})^{-2\alpha_{k}}\left(\left\Vert [u^{k}]\right\Vert _{_{0,\tilde{\gamma_{s}}}}^{2}+\left\Vert [(u_{\tau_{k}}^{k})]\right\Vert _{_{1/2,\tilde{\gamma_{s}}}}^{2}+\left\Vert [(u_{\theta_{k}}^{k})]\right\Vert _{_{1/2,\tilde{\gamma_{s}}}}^{2}\right)\nonumber \\
 & + & \sum_{k=0}^{4}\sum_{\gamma_{s}\subseteq L_{k}}\!\!\!d(E_{k},\gamma_{s})^{-2\alpha_{k}}\left(\left\Vert [u^{k}]\right\Vert _{_{\frac{3}{2},\tilde{\gamma}_{s}}}^{2}+\left\Vert \left[p\frac{\partial u^{k}}{\partial\theta_{k}}\right]\right\Vert _{\frac{1}{2},\tilde{\gamma}_{s}}^{2}\right)\nonumber \\
 & + & \sum_{m\in\mathcal{D}}\sum_{k=m-1}^{m}\sum_{{\gamma_{s}\subseteq\partial\varPi^{k}\cap\Gamma_{m},\atop \mu(\tilde{\gamma_{s}})<\infty}}\!\!\!\!d(E_{k},\gamma_{s})^{-2\alpha_{k}}\left(\left\Vert (u^{k}-h_{k})-(l_{m-k+1}^{m}-a_{k})\right\Vert _{_{0,\tilde{\gamma_{s}}}}^{2}\right.\label{eq:6}\\
 & + & \left.\left\Vert u_{\tau_{k}}^{k}-(l_{m-k+1}^{m})_{\tau_{k}}\right\Vert _{_{1/2,\tilde{\gamma_{s}}}}^{2}\right)+\sum_{m\in\mathcal{D}}\sum_{k=m-1}^{m}(h_{k}-a_{k})^{2}\nonumber \\
 & + & \sum_{m\in\mathcal{N}}\sum_{k=m-1}^{m}\sum_{{\gamma_{s}\subseteq\partial\varPi^{k}\cap\Gamma_{m},\atop \mu(\tilde{\gamma_{s}})<\infty}}\!\!\!d(E_{k},\gamma_{s})^{-2\alpha_{k}}\left\Vert \left(\frac{\partial u^{k}}{\partial n}\right)-l_{m-k+1}^{m}\right\Vert _{_{1/2,\tilde{\gamma_{s}}}}^{2}.\nonumber 
\end{eqnarray}

In the above $\mu(\tilde{\gamma_{s}})$ denotes the measure of $\tilde{\gamma_{s}}.$ 

Define
\begin{eqnarray}
\mathfrak{\mathcal{\mathfrak{r}}}_{interior}^{^{N,W}}\left(\mathcal{F}_{u}\right) & = & \sum_{k=0}^{4}\sum_{j=N+1}^{J_{k}}\sum_{i=1}^{I_{k}}\left\Vert (\mathcal{L}_{i,j}^{k})u_{i,j}^{k}(\xi,\eta)-F_{i,j}^{k}\left(\xi,\eta\right)\right\Vert _{_{0,\Omega_{i,j}^{k}}}^{2}\nonumber \\
 & + & \sum_{\gamma_{s}\subseteq\Omega^{0}\cup\Omega^{1},\gamma_{s}\varsubsetneq L_{k}}\left(\left\Vert [u^{k}]\right\Vert _{_{0,\gamma_{s}}}^{2}+\left\Vert [(u_{x_{1}}^{k})]\right\Vert _{_{1/2,\gamma_{s}}}^{2}+\left\Vert [(u_{x_{2}}^{k})]\right\Vert _{_{1/2,\gamma_{s}}}^{2}\right)\label{eq:-3}\\
 & + & \sum_{\gamma_{s}\subseteq L_{k}\subseteq\Omega^{0}\cup\Omega^{1}}\left(\left\Vert [u^{k}]\right\Vert _{_{\frac{3}{2},\gamma_{s}}}^{2}+\left\Vert \left[\left(p\frac{\partial u^{k}}{\partial n}\right)\right]\right\Vert _{_{1/2,\gamma_{s}}}^{2}\right)\nonumber \\
 & + & \sum_{l\in\mathcal{D}}\sum_{\gamma_{s}\subseteq\partial\Omega^{1}\cap\Gamma_{l}}\left(\left\Vert u^{k}-o^{l,k}\right\Vert _{_{0,\gamma_{s}}}^{2}+\left\Vert \left(\frac{\partial u^{k}}{\partial T}\right)-\left(\frac{\partial o^{l,k}}{\partial T}\right)\right\Vert _{_{1/2,\gamma_{s}}}^{2}\right)\nonumber \\
 & + & \sum_{l\in\mathcal{N}}\sum_{\gamma_{s}\subseteq\partial\Omega^{1}\cap\Gamma_{l}}\left\Vert \left(\frac{\partial u^{k}}{\partial n}\right)-o^{l,k}\right\Vert _{_{1/2,\gamma_{s}}}^{2}.\nonumber 
\end{eqnarray}
 Let
\begin{eqnarray*}
 &  & \mathfrak{\mathfrak{\mathcal{\mathfrak{r}}}}^{^{N,W}}(\mathcal{F}_{u})=\mathfrak{\mathcal{\mathfrak{r}}}_{vertices}^{^{N,W}}(\mathcal{F}_{u})+\mathcal{\mathfrak{r}}_{interior}^{^{N,W}}(\mathcal{F}_{u}).
\end{eqnarray*}

We choose as our approximate solution the unique $\mathcal{F}_{z}\in\Pi^{N,W},$
the space of spectral element functions, which minimizes the functional
$\mathfrak{\mathcal{\mathcal{\mathfrak{r}}}}^{^{N,W}}(\mathcal{F}_{u})$
over all $\mathcal{F}_{u}.$

The method is essentially a least-squares method and the solution
is obtained at Gauss-Legendre-Lobatto points using preconditioned
conjugate gradient method without storing the stiffness matrix and
load vector. The residuals in the normal equations can be computed
efficiently and inexpensively as shown in \cite{KISHORENAGA,tomar}. 

The minimization leads to the normal equations 
\begin{eqnarray*}
 &  & AU=h.
\end{eqnarray*}
The vector $U$ composed of the values of the spectral element functions
at Gauss-Legendre-Lobatto points is divided into two sub vectors one
consisting of the values of the spectral element functions at the
vertices of the domain constitute the set of common boundary values
$U_{B},$ and the other consisting of the remaining values which we
denote by $U_{I}.$ The computation of $U_{I}$ and $U_{B}$ is described
in \cite{KISHORENAGA,tomar}. 

An efficient preconditioner has been used which is proposed in \cite{pravirpankaj}
for the matrix $A$ so that the condition number of the preconditioned
system is as small as possible. The condition number of the preconditioned
system is $O((ln\,W)^{2})$. The preconditioner is a block diagonal
matrix, where each diagonal block is constructed using the separation
of variable technique. So the solution is obtained to an exponential
accuracy using $O(W\,ln\,W)$ iterations of the PCGM. After obtaining
the nonconforming solution at the Gauss-Legendre-Lobatto points, a
set of corrections are performed \cite{schawb} so that the solution
is conforming and belongs to $H^{1}(\Omega).$ 

Then for $W$ large enough the error estimate
\begin{eqnarray*}
 &  & \left\Vert u-u_{ap}\right\Vert _{1,\Omega}\leq C\,e^{-bW}
\end{eqnarray*}
holds, where $C$ and $b$ are constants and $u_{ap}$ is the corrected
solution.

\section{Numerical Results}

Here we consider few numerical examples to show the effectiveness
of the proposed method. For simplicity, we have considered $W_{j}=W$
for all $j$ and the number of layers $N=W$ in the geometric mesh.
The relative error $\left\Vert e\right\Vert _{ER}=\frac{\left\Vert e\right\Vert _{1}}{\left\Vert u\right\Vert _{1}},$
where $e=u-u_{ap}$ is the difference in the exact solution $u$ and
the approximate solution $u_{ap}$ measured in $H^{1}$ norm. ``Iters''
is the total number of iterations to compute $U_{I}$ and $U_{B}.$\\
\\
\textbf{Example 1: Interface problem with singularity at the intersection
of an interface and the boundary}\\
\\
Let us consider the interface problem on the domain $\Omega=\Omega_{1}\cup\Omega_{2}$
as shown in Fig. 6
\begin{eqnarray*}
 &  & -\nabla.(p\nabla u)=0\,\,\textrm{in}\,\,\Omega\\
 &  & u=0\,\,\textrm{on}\,\,\Gamma_{1}\\
 &  & \frac{\partial u}{\partial\theta}=0\,\,\textrm{on}\,\,\Gamma_{2}
\end{eqnarray*}
 where the coefficient $p$ is piecewise constant:
\begin{eqnarray*}
 &  & p=\begin{cases}
1 & \textrm{in}\,\,\Omega_{1}\\
p & \textrm{in}\,\,\Omega_{2}.
\end{cases}
\end{eqnarray*}

\begin{figure}[H]
~~~~~~~~~~~~~~~~~~~~~~~~~~~~~~~~~~~~~~~~~~~~~~~~~~~~~~~~~~~~~~~~~~~~~~~~\includegraphics{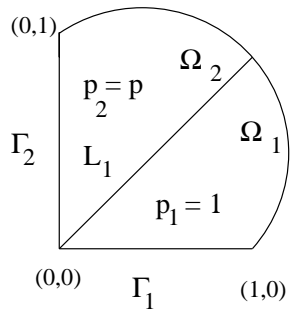}

\caption{The interface problem on a sector}
\end{figure}

Let $(r,\theta)$ be the polar coordinates centered at the origin
$(0,0).$ Assume that the interface conditions are satisfied at $\theta=\frac{\pi}{4}$.
That is, on $L_{1}$
\begin{eqnarray*}
 &  & u(r,\frac{\pi}{4}-0)=u(r,\frac{\pi}{4}+0)\\
 &  & \frac{\partial u}{\partial\theta}(r,\frac{\pi}{4}-0)=p\frac{\partial u}{\partial\theta}(r,\frac{\pi}{4}+0).
\end{eqnarray*}

Let the solution of the above interface problem be in the form $r^{\lambda}W(\theta).$
As explained in \cite{guooh,ohbabu}, $\lambda$ and $W(\theta)$
can be obtained by solving the following Sturm-Liouville problem corresponding
to the above interface problem:
\begin{eqnarray}
 &  & \frac{d}{d\theta}(p\frac{dW}{d\theta})+\lambda^{2}pW=0\,\,\textrm{in}\,\,\Omega\label{eq:-1}\\
 &  & W(0)=0,\frac{dW}{d\theta}(\frac{\pi}{2})=0.\nonumber 
\end{eqnarray}
 The function $W(\theta)$ required to satisfy 
\begin{eqnarray}
 &  & W(\frac{\pi}{4}-0)=W(\frac{\pi}{4}+0)\label{eq:-2}\\
 &  & \frac{dW}{d\theta}(\frac{\pi}{4}-0)=p\frac{dW}{d\theta}(\frac{\pi}{4}+0).\nonumber 
\end{eqnarray}
 The solution of the above differential equation $W$ is of the form
\begin{eqnarray*}
 &  & W(\theta)=\begin{cases}
 & C_{1}\,\textrm{cos}\lambda\theta+C_{2}\,\textrm{sin}\lambda\theta\,\,\,\textrm{in}\,\,\Omega_{1}\\
 & C_{3}\,\textrm{cos}\lambda\theta+C_{4}\,\textrm{sin}\lambda\theta\,\,\,\textrm{in}\,\,\Omega_{2}
\end{cases}
\end{eqnarray*}
Therefore the solution of the interface problem has the following
form
\begin{eqnarray*}
 &  & u_{1}=r^{\lambda}(C_{1}\,\textrm{cos}\lambda\theta+C_{2}\,\textrm{sin}\lambda\theta)\,\,\,\textrm{in}\,\,\Omega_{1}\\
 &  & u_{2}=r^{\lambda}(C_{3}\,\textrm{cos}\lambda\theta+C_{4}\,\textrm{sin}\lambda\theta)\,\,\,\textrm{in}\,\,\Omega_{2}.
\end{eqnarray*}
 The constants $C_{1},C_{2},C_{3}\,\,\textrm{and}\,\,C_{4}$ and the
eigenvalues $\lambda$ can be obtained using the above defined boundary
and interface conditions (8) and (9).

Now, after applying the boundary conditions, we get 
\begin{eqnarray*}
 &  & W(0)=0\rightarrow C_{1}=0\\
 &  & \frac{dW}{d\theta}(\frac{\pi}{2})=0\rightarrow-C_{3}\,\textrm{sin}\frac{\lambda\pi}{2}+C_{4}\,\textrm{cos}\frac{\lambda\pi}{2}=0.
\end{eqnarray*}
 The interface conditions gives 
\begin{eqnarray*}
 &  & W(\frac{\pi}{4}-)=W(\frac{\pi}{4}+)\\
 &  & \Rightarrow C_{2}\,\textrm{sin}\frac{\lambda\pi}{4}-C_{3}\,\textrm{cos}\frac{\lambda\pi}{4}-C_{4}\,\textrm{sin}\frac{\lambda\pi}{4}=0
\end{eqnarray*}
and 
\begin{eqnarray*}
 &  & \frac{dW}{d\theta}(\frac{\pi}{4}-0)=p\frac{dW}{d\theta}(\frac{\pi}{4}+0)\\
 &  & \Rightarrow C_{2}\,\textrm{cos}\frac{\lambda\pi}{4}+C_{3}\,p\textrm{sin}\frac{\lambda\pi}{4}-C_{4}\,p\textrm{cos}\frac{\lambda\pi}{4}=0.
\end{eqnarray*}
 So we have the following homogeneous system 
\begin{eqnarray}
 &  & \left[\begin{array}{ccc}
0 & -\textrm{sin}\frac{\lambda\pi}{2} & \textrm{cos}\frac{\lambda\pi}{2}\\
\textrm{sin}\frac{\lambda\pi}{4} & -\textrm{cos}\frac{\lambda\pi}{4} & -\textrm{sin}\frac{\lambda\pi}{4}\\
\textrm{cos}\frac{\lambda\pi}{4} & p\textrm{sin}\frac{\lambda\pi}{4} & -p\textrm{cos}\frac{\lambda\pi}{4}
\end{array}\right]\left[\begin{array}{c}
C_{2}\\
C_{3}\\
C_{4}
\end{array}\right]=\left[\begin{array}{c}
0\\
0\\
0
\end{array}\right].\label{eq:}
\end{eqnarray}
 In order for the system of unknowns $C_{2},C_{3},C_{4}$ to have
a non-trivial solution, the determinant of the coefficient matrix
$A$ of the system should be zero.
\begin{eqnarray*}
 &  & \left|\begin{array}{ccc}
0 & -\textrm{sin}\frac{\lambda\pi}{2} & \textrm{cos}\frac{\lambda\pi}{2}\\
\textrm{sin}\frac{\lambda\pi}{4} & -\textrm{cos}\frac{\lambda\pi}{4} & -\textrm{sin}\frac{\lambda\pi}{4}\\
\textrm{cos}\frac{\lambda\pi}{4} & p\textrm{sin}\frac{\lambda\pi}{4} & -p\textrm{cos}\frac{\lambda\pi}{4}
\end{array}\right|=0\\
 &  & \Longrightarrow\frac{(1-p)}{2}\textrm{si\ensuremath{n^{2}}}\frac{\lambda\pi}{2}+\textrm{cos}\frac{\lambda\pi}{2}+(p-1)\textrm{si\ensuremath{n^{2}}}\frac{\lambda\pi}{4}cos\frac{\lambda\pi}{2}=0\\
 &  & \Longrightarrow\frac{(1-p)}{2}\textrm{si\ensuremath{n^{2}}}\frac{\lambda\pi}{2}+\textrm{cos}\frac{\lambda\pi}{2}+\frac{(p-1)}{2}(1-\textrm{cos}\frac{\lambda\pi}{2})\textrm{cos}\frac{\lambda\pi}{2}=0\\
 &  & \Longrightarrow\frac{(1-p)}{2}+\left[\frac{(p-1)}{2}+1\right]\textrm{cos}\frac{\lambda\pi}{2}=0.
\end{eqnarray*}

By solving $\frac{(1-p)}{2}+\left[\frac{(p-1)}{2}+1\right]\textrm{cos}\frac{\lambda\pi}{2}=0,$
we obtain the eigenvalues $\lambda_{k}$ which are positive real values.
The smallest eigenvalue $\lambda_{0}$ among $\lambda_{k}$ gives
the value of the exponent in the leading order singular term in the
expansion of $u.$ The value of $\lambda_{0}$ for different values
of $p$ is given in the following Table 1. The strength of the singularity
increases as $p$ increases. 

\begin{table}[H]
~~~~~~~~~~~~~~~~~~~~~~~~~~~~~~~~~~~~~~~~~~~~~~~~~~~~~~~~~~%
\begin{tabular}{cc}
\hline 
$p$ & $\lambda_{0}$\tabularnewline
\hline 
5 & 0.53544092\tabularnewline
10 & 0.38996444\tabularnewline
30 & 0.22992823\tabularnewline
50 & 0.1788770\tabularnewline
100 & 0.12690206\tabularnewline
\hline 
\end{tabular}

\caption{The exponent of leading order singular term for different $p$ }
\end{table}

Now, we find the constants $C_{2},C_{3}$ and $C_{4}.$ From the above
linear system 
\begin{eqnarray*}
 &  & C_{4}=C_{3}\frac{\textrm{sin}\frac{\lambda\pi}{2}}{\textrm{cos}\frac{\lambda\pi}{2}}=C_{3}\,\textrm{tan}\frac{\lambda\pi}{2}.
\end{eqnarray*}
 We choose $C_{3}=1.$ Therefore $C_{4}=tan\frac{\lambda\pi}{2}.$
Then, one can easily find the value of $C_{2}$ from any one of the
equations in the linear system. The value of $C_{2}$ is given by
\begin{eqnarray*}
 &  & C_{2}=\textrm{cot}\frac{\lambda\pi}{4}+\textrm{tan}\frac{\lambda\pi}{2}.
\end{eqnarray*}
 Therefore the leading order singular term in the expansion of $u$
has the following form
\begin{eqnarray*}
 &  & u_{1}=r^{\lambda_{0}}(\textrm{cot}\frac{\lambda_{0}\pi}{4}+\textrm{tan}\frac{\lambda_{0}\pi}{2})\textrm{sin}\lambda_{0}\theta\,\,\,\textrm{in}\,\,\Omega_{1}\\
 &  & u_{2}=r^{\lambda_{0}}(\textrm{cos}\lambda_{0}\theta+\textrm{tan}\frac{\lambda_{0}\pi}{2}\textrm{sin}\lambda_{0}\theta)\,\,\,\textrm{in}\,\,\Omega_{2}.
\end{eqnarray*}
 \textbf{\textit{Remark:}}\textit{ The solution $u$ has singular
behavior at the point $(0,0)$ and the strength of the singularity
is very strong for larger values of $p$. These singularities are
more stronger than the singularities which generally arises in elliptic
problems due to the nonsmooth domains. \medskip{}
}

Now we present the numerical solution of this problem. We consider
the Dirichlet boundary condition on $\rho=1$ ( Fig. 6). Since the
strength of the singularity is very strong at the corners a very refined
mesh as well as higher degree of approximation is needed to get a
good accuracy. In \cite{guooh} the numerical solution is obtained
using $hp$ finite element method. They have used a geometric mesh
near the corner with geometric ratio $0.15$ and tabulated the relative
error for different values of the degree of approximation $W$ with
$2W$ layers in the geometric mesh in the radial direction. 

In the following table we have presented the numerical results for
$p=5.$ As explained above the exponent in the leading order singular
term in the solution for $p=5$ is 0.53544092. We consider the geometric
ratio $\mu=0.15.$ The relative error is obtained for different values
of $W$ and shown in the following Table 2. Table 2 also shows the
number of iterations. 

\begin{table}[H]
~~~~~~~~~~~~~~~~~~~~~~~~~~~~~~~~~~~~~~~~~~~~~~~~~~~~~~~
~~~~~~~~%
\begin{tabular}{ccc}
\hline 
$W$ & $\left\Vert e\right\Vert _{ER}\%$ & Iters\tabularnewline
\hline 
2 & 11.254417 & 37\tabularnewline
3 & 4.29850 & 77\tabularnewline
4 & 1.541124 & 118\tabularnewline
5 & 0.5575801 & 159\tabularnewline
6 & 0.2017785 & 204\tabularnewline
7 & 0.0730631 & 250\tabularnewline
8 & 0.0264555 & 291\tabularnewline
9 & 0.0095797 & 335\tabularnewline
\hline 
\end{tabular}

\caption{The relative error and iterations against $W$}
\end{table}

Fig. 7 shows the log of the relative error against the degree of approximation
$W$ and the relation is linear. This shows the exponential accuracy
of the method. 
\begin{figure}[H]
\includegraphics[scale=0.5]{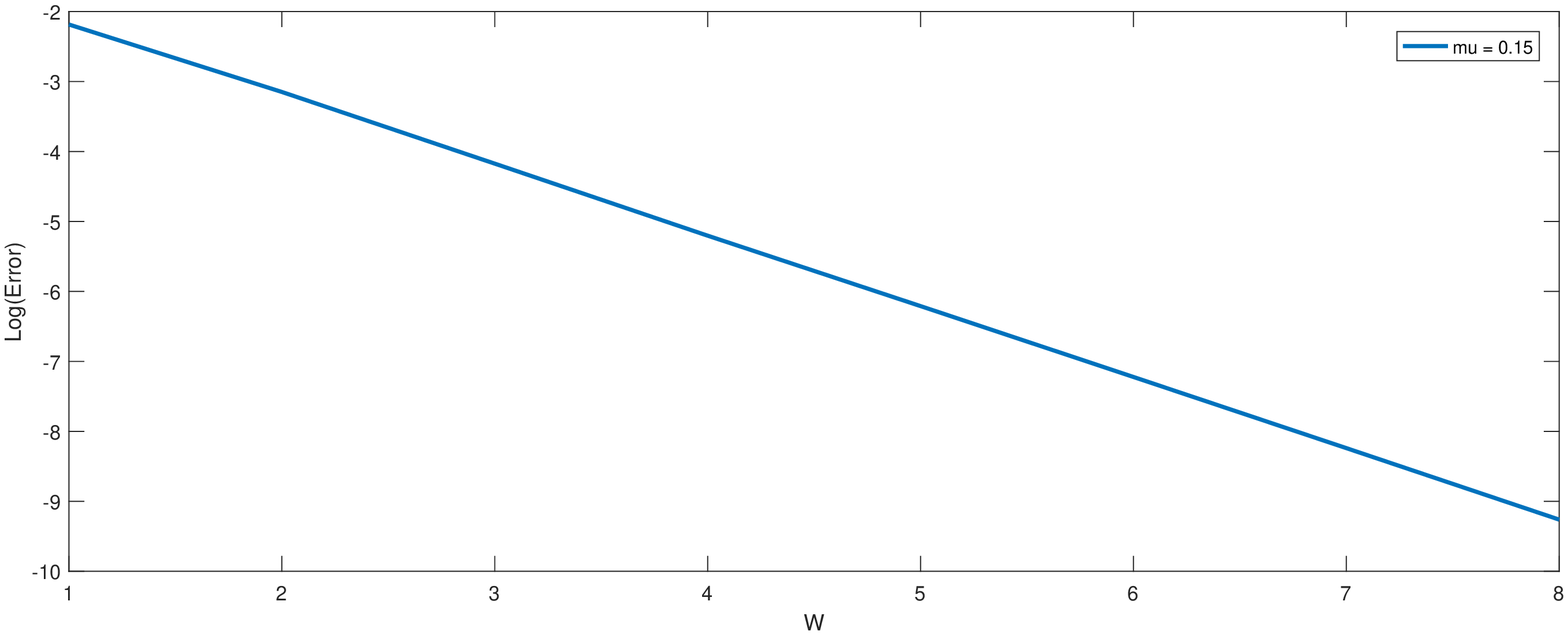}

\caption{Log of $\left\Vert e\right\Vert _{ER}$ against $W$}

\end{figure}

Now we consider $p=10.$ The exponent in the leading order singular
term of the solution is 0.38996444. This is strong compared to the
previous exponent. So we need more refined grid near the singular
point. Here we consider the geometric ratio $\mu=0.15$ and $\mu=e^{-\pi}.$
The relative error and iteration count for different values of $W$
is tabulated in the Table 3.

\begin{table}[H]
~~~~~~~~~~~~~~~~~~~~~~~~~~~~~~~~~~~~~~~%
\begin{tabular}{ccc|cc}
\hline 
 & $\mu=0.15$ &  &  & $\mu=e^{-\pi}$\tabularnewline
\hline 
$W$ & $\left\Vert e\right\Vert _{ER}\%$ & Iters & $\left\Vert e\right\Vert _{ER}\%$ & Iters\tabularnewline
\hline 
2 & 19.08735 & 44 & 7.19601 & 48\tabularnewline
3 & 9.632562 & 95 & 2.300786 & 122\tabularnewline
4 & 4.66738 & 159 & 0.676060 & 208\tabularnewline
5 & 2.24022 & 226 & 0.198810 & 279\tabularnewline
6 & 1.07102 & 290 & 0.058287 & 368\tabularnewline
7 & 0.511375 & 346 & 0.017128 & 449\tabularnewline
8 & 0.244061 & 424 & 0.00503136 & 525\tabularnewline
9 & 0.116472 & 474 & 0.00147676 & 623\tabularnewline
\hline 
\end{tabular}

\caption{The relative error and iterations against $W$ for $\mu=0.15$ and
$\mu=e^{-\pi}$}
\end{table}

The error decays slowly for the geometric ratio $\mu=0.15.$ One can
get better accuracy by increasing the number of layers in the geometric
mesh. But this increases the number of degrees of freedom. For $\mu=e^{-\pi}$
the error decays very fast with an increase in the iteration count.
Even better accuracy can be achieved with the geometric ratio $\mu=e^{-1.5\pi}.$
In the Figure 8 the graph of log of relative error vs. $W$ has been
drawn for $\mu=0.15$ and $\mu=e^{-\pi}.$ The relation is linear. 

\begin{figure}[H]
\includegraphics[scale=0.5]{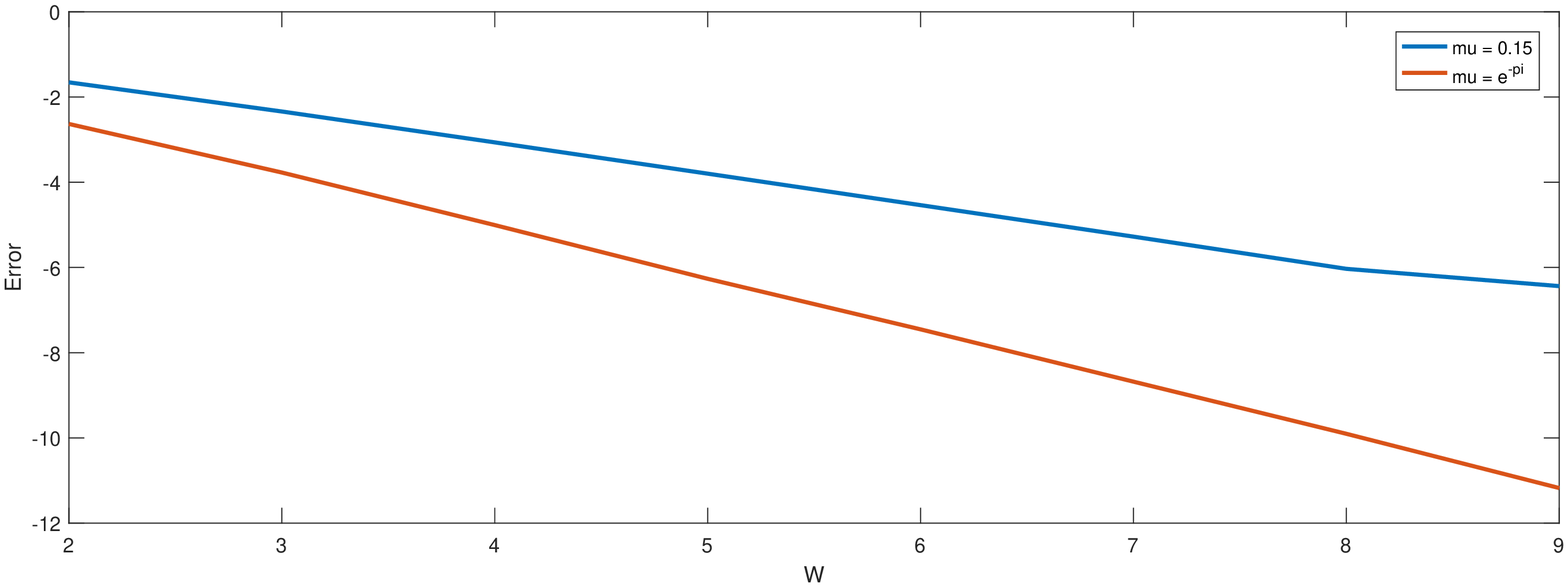}

\caption{Log of $\left\Vert e\right\Vert _{ER}$ against $W$}

\end{figure}

Now we consider $p=30.$ In this case the exponent in the leading
order singular term of the solution is 0.22992823. Here we consider
four different geometric mesh with ratio $\mu=0.15,$$\mu=0.15$ with
more number of layers (just double of the degree of the approximation
$W$) in radial direction, $\mu=e^{-\pi}$ and $\mu=e^{-1.5\pi}.$
The relative error and iterations are shown for different values of
$W$ and for different geometric ratios in Table 4. 

\begin{table}[H]
~~~~~~~~~~~~~~~~%
\begin{tabular}{ccc|cc|cc|cc}
\hline 
 & $\mu=0.15$ &  & $\mu=0.15$ &  & $\mu=e^{-\pi}$ &  & $\mu=e^{-1.5\pi}$ & \tabularnewline
\hline 
$W$ & $\left\Vert e\right\Vert _{ER}\%$ & Iters & $\left\Vert e\right\Vert _{ER}\%$ & Iters & $\left\Vert e\right\Vert _{ER}\%$ & Iters & $\left\Vert e\right\Vert _{ER}\%$ & Iters\tabularnewline
\hline 
2 & 32.48562 & 49 & 14.36468 & 103 & 18.848278 & 60 & 9.12371 & 63\tabularnewline
3 & 22.58167 & 122 & 6.68331 & 244 & 10.28457 & 156 & 3.53691 & 169\tabularnewline
4 & 15.27828 & 212 & 2.83311 & 402 & 5.123292 & 272 & 1.210048 & 326\tabularnewline
5 & 10.16691 & 308 & 1.118838 & 551 & 2.50895 & 404 & 0.409445 & 500\tabularnewline
6 & 6.68372 & 409 & 0.496950 & 721 & 1.22152 & 548 & 0.138267 & 713\tabularnewline
7 & 4.36046 & 502 & 0.207736 & 891 & 0.593652 & 667 & 0.046691 & 898\tabularnewline
8 &  &  &  &  &  &  & 0.015766 & 1051\tabularnewline
9 &  &  &  &  &  &  & 0.005324 & 1263\tabularnewline
\hline 
\end{tabular}

\caption{The relative error and iterations against $W$ }
\end{table}

The results shows the geometric ratio $\mu=e^{-1.5\pi}$ gives better
results. Even better accuracy can be achieved with the geometric ratio
$\mu=e^{-2\pi}$ with an increase in the number of iterations. The
Fig. 9 shows the graph of log of relative error against $W$ for different
values of $\mu$. The relation is linear in all cases but the convergence
is faster for $\mu=e^{-1.5\pi}$. 

\begin{figure}[H]
\includegraphics[scale=0.5]{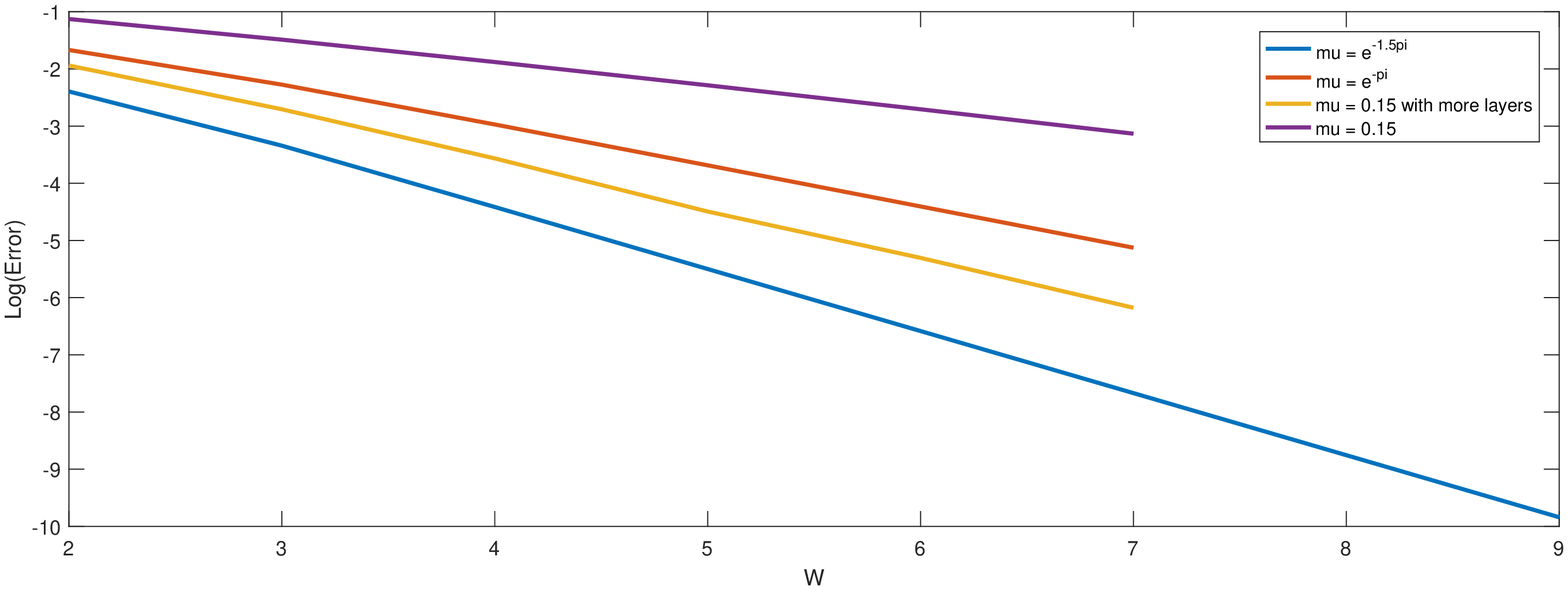}

\caption{Log of $\left\Vert e\right\Vert _{ER}$ against $W$}

\end{figure}

Now consider $p=50,100.$ The exponents in the leading order singular
term of the solution are 0.1788770 and 0.12690206 respectively. So
we need very refined mesh in the neighbourhood of the singular point.
So we consider the geometric ratio $\mu=e^{-2\pi}.$ The relative
error and iterations are tabulated for different values of $W$ in
Table 5. The numerical results shows the good performance of the method.\\
\begin{table}[H]
~~~~~~~~~~~~~~~~~~~~~~~~~~~~~~~~~~~~~~~~~~~~~~~~~~%
\begin{tabular}{ccc|cc}
\hline 
 & $p=50$ &  &  & $p=100$\tabularnewline
\hline 
$W$ & $\left\Vert e\right\Vert _{ER}\%$ & Iters & $\left\Vert e\right\Vert _{ER}\%$ & Iters\tabularnewline
\hline 
2 & 8.353101 & 80 & 15.9732 & 92\tabularnewline
3 & 3.119198 & 219 & 8.29721 & 241\tabularnewline
4 & 1.019195 & 438 & 3.83097 & 475\tabularnewline
5 & 0.330973 & 700 & 1.73125 & 778\tabularnewline
6 & 0.106907 & 1030 & 0.777691 & 1137\tabularnewline
7 & 0.034519 & 1354 & 0.348880 & 1486\tabularnewline
8 & 0.011145 & 1631 & 0.156460 & 1806\tabularnewline
9 & 0.003598 & 1967 & 0.0701626 & 2358\tabularnewline
\hline 
\end{tabular}

\caption{The relative error and iterations against $W$}
\end{table}

Fig. 10 shows the graph of log relative error against $W$ for $p=50$
and $p=100.$ The relation is linear. This show the exponential convergence
of the proposed method. 

\begin{figure}[H]
\includegraphics[scale=0.5]{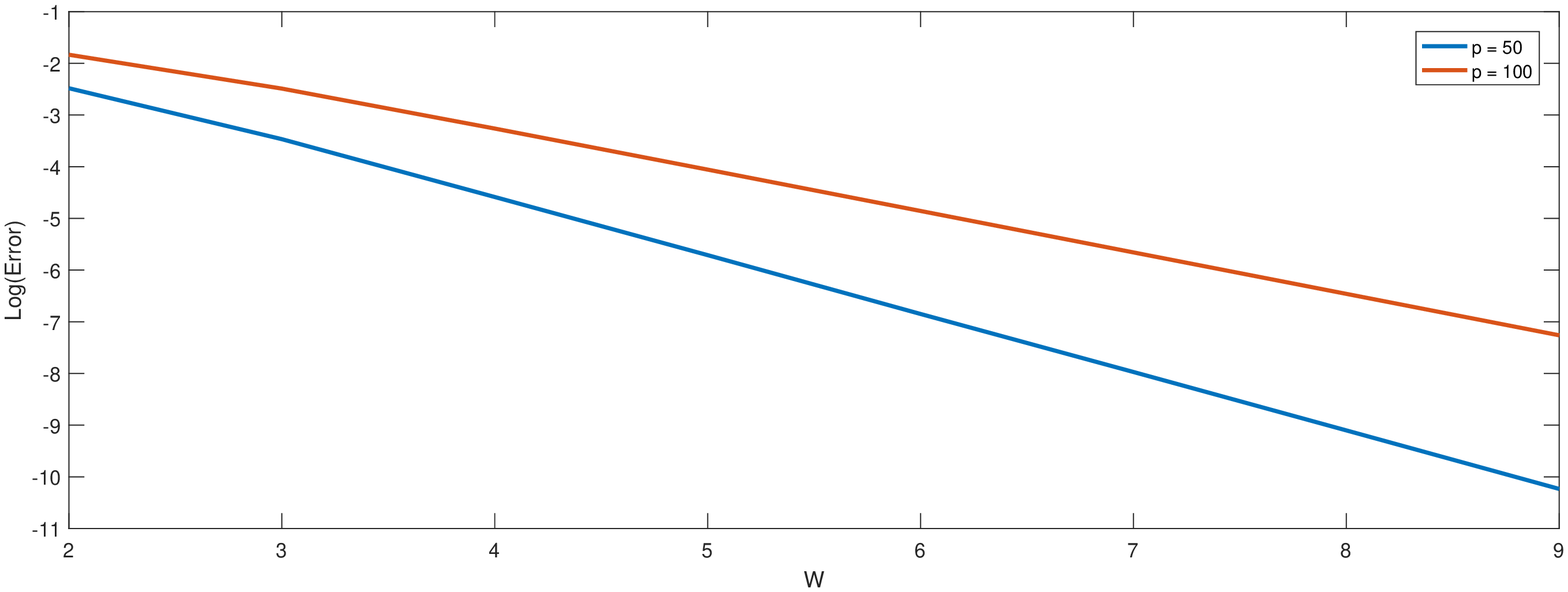}

\caption{Log of $\left\Vert e\right\Vert _{ER}$ against $W$}

\end{figure}

\bigskip{}
\bigskip{}
\textbf{$\!\!\!\!\!\!\!\!$Example 2: Interface problem with singularity
at the intersection of two interfaces}\bigskip{}

\begin{figure}[H]
~~~~~~~~~~~~~~~~~~~~~~~~~~~~~~~~~~~~~~~~~~~~~~~~~~~~~~~~~~~~~~~~~~~~~~\includegraphics{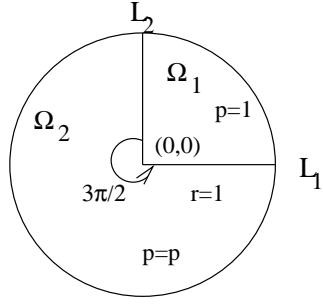}

\caption{The domain $\Omega$ with $L_{1}$and $L_{2}$ meet at $(0,0)$}
\end{figure}

Consider the following interface problem on the domain $\Omega$ as
shown in Fig. 11. 

\begin{eqnarray*}
 &  & -\nabla.(p\nabla u)=0\,\,\textrm{in}\,\,\Omega
\end{eqnarray*}
 where the coefficient $p$ is piecewise constant:
\begin{eqnarray*}
 &  & p=\begin{cases}
1 & \textrm{in}\,\,\Omega_{1}\\
p & \textrm{in}\,\,\Omega_{2}
\end{cases}
\end{eqnarray*}
with Dirichlet boundary condition on the circle of radius $1.$ Assume
that the two interfaces $L_{1}=\{(r,0),0\leq r\leq1\}$ and $L_{2}=\{(r,\frac{\pi}{2}),0\leq r\leq1\}$
meets at the point $E_{0}=(0,0)$ and $u$ satisfies the interface
conditions 
\begin{eqnarray*}
 &  & u(\theta=0)=u(\theta=2\pi)\,\,\textrm{and}\,\,\frac{\partial u}{\partial\theta}(0)=p\frac{\partial u}{\partial\theta}(2\pi)\\
 &  & u(\frac{\pi}{2}-)=u(\frac{\pi}{2}+)\,\,\textrm{and}\,\,\frac{\partial u}{\partial\theta}(\frac{\pi}{2}-)=p\frac{\partial W}{\partial\theta}(\frac{\pi}{2}+).
\end{eqnarray*}

Here we are only interested in the behavior of the solution at $(0,0).$
So as explained in Section 2, we need to solve the Sturm-Liouville
problem
\[
\frac{d}{d\theta}(p\frac{dW}{d\theta})+\lambda^{2}pW=0\,\,\textrm{in}\,\,\Omega
\]
with 
\begin{eqnarray*}
 &  & W(\theta=0)=W(\theta=2\pi)\,\textrm{\,and}\,\,\frac{dW}{d\theta}(0)=p\frac{dW}{d\theta}(2\pi)\\
 &  & W(\frac{\pi}{2}-)=W(\frac{\pi}{2}+)\,\,\textrm{and}\,\,\frac{dW}{d\theta}(\frac{\pi}{2}-)=p\frac{dW}{d\theta}(\frac{\pi}{2}+).
\end{eqnarray*}
 The solution of the above differential equation $W$ is of the form
\begin{eqnarray*}
 &  & W(\theta)=\begin{cases}
 & C_{1}\,\textrm{cos}\lambda\theta+C_{2}\,\textrm{sin}\lambda\theta\,\,\,\textrm{in}\,\,\Omega_{1}\\
 & C_{3}\,\textrm{cos}\lambda\theta+C_{4}\,\textrm{sin}\lambda\theta\,\,\textrm{\,in}\,\,\Omega_{2}.
\end{cases}
\end{eqnarray*}
 As explained in the above example, we get a homogeneous system of
equations in unknowns $C_{1},C_{2},C_{3}\,\,\textrm{and}\,\,C_{4}.$
In order to have a non-trivial solution, the determinant of the coefficient
matrix $A$ of the system should be zero. This gives an equation in
$\lambda$ and the eigenvalues $\lambda_{k}^{'}$s are the solutions
of this equation. We have obtained the smallest eigenvalue $\lambda_{0}$
for different values of $p$ and tabulated in the following Table
6. 

\begin{table}[H]
~~~~~~~~~~~~~~~~~~~~~~~~~~~~~~~~~~~~~~~~~~~~~~~~~~%
\begin{tabular}{cc}
\hline 
$p$ & $\lambda_{0}$\tabularnewline
\hline 
5 & 0.783653104062978\tabularnewline
10 & 0.731691778699314\tabularnewline
30 & 0.690135330693010\tabularnewline
50 & 0.680988694144617\tabularnewline
100 & 0.673921228717518\tabularnewline
500 & 0.668132968861755\tabularnewline
\hline 
\end{tabular}

\caption{The exponent of leading order singular term for different $p$}
\end{table}

The singularities in this case are not so strong as the singularities
which we have seen in example 1. 

We obtain the constants $C_{1},C_{2},C_{3}$ and $C_{4}$ using the
above interface conditions. 
\begin{eqnarray*}
 &  & W(\theta=0)=W(\theta=2\pi)\\
 &  & \Longrightarrow C_{1}=C_{3}\textrm{cos}2\pi\lambda+C_{4}\textrm{sin}2\pi\lambda\\
 &  & \frac{dW}{d\theta}(0)=p\frac{dW}{d\theta}(2\pi)\\
 &  & \Longrightarrow C_{2}=-pC_{3}\textrm{sin}2\pi\lambda+pC_{4}\textrm{cos}2\pi\lambda.
\end{eqnarray*}
 Now let $C_{4}=1.$ Then $C_{1}=C_{3}\textrm{cos}2\pi\lambda+\textrm{sin}2\pi\lambda$
and $C_{2}=-C_{3}p\textrm{sin}2\pi\lambda+p\textrm{cos}2\pi\lambda.$
Then one can easily find $C_{3}.$ The value of $C_{3}$ is given
by
\[
C_{3}=\frac{(\textrm{sin}\frac{\lambda\pi}{2}-p\,\textrm{cos}2\pi\lambda\,\textrm{sin}\frac{\lambda\pi}{2}-\textrm{sin}2\pi\lambda\,\textrm{cos}\frac{\lambda\pi}{2})}{(\textrm{cos}2\pi\lambda\textrm{\,cos}\frac{\lambda\pi}{2}-p\textrm{sin}2\pi\lambda\,\textrm{sin}\frac{\lambda\pi}{2}-\textrm{cos}\frac{\lambda\pi}{2})}.
\]

Therefore the leading order singular term in the solution of the interface
problem is given by 
\begin{eqnarray*}
 &  & u_{1}=r^{\lambda_{0}}(C_{1}\textrm{cos}\lambda_{0}\theta+C_{2}\textrm{sin}\lambda_{0}\theta)\,\,\,\textrm{in}\,\,\Omega_{1}\\
 &  & u_{2}=r^{\lambda_{0}}(C_{3}\textrm{cos}\lambda_{0}\theta+\textrm{sin}\lambda_{0}\theta)\,\,\,\textrm{in}\,\,\Omega_{2}
\end{eqnarray*}
 with the constants $C_{1},C_{2}$ and $C_{3}$ given above.

We have obtained the numerical solution for $p=500.$ Table 7 shows
the relative error and iterations for different values of $W.$ 

\begin{table}[H]
~~~~~~~~~~~~~~~~~~~~~~~~~~~~~~~~~~~~~~~~~~~~~~~~~~~~~~~~~~~~~~%
\begin{tabular}{ccc}
\hline 
$W$ & $\left\Vert e\right\Vert _{ER}\%$ & Iters\tabularnewline
\hline 
2 & 28.5515688 & 42\tabularnewline
3 & 2.16885204 & 170\tabularnewline
4 & 0.58543007 & 306\tabularnewline
5 & 0.16244907 & 467\tabularnewline
6 & 0.04476168 & 674\tabularnewline
7 & 0.01260480 & 840\tabularnewline
8 & 0.00354783 & 997\tabularnewline
9 & 0.00099866 & 1215\tabularnewline
\hline 
\end{tabular}

\caption{The relative error for different values of $W$}

\end{table}

Figure 12 shows the graph of log of relative error against $W$ for
$p=500.$ The relation is linear. This shows the exponential accuracy
of the method.

\begin{figure}[H]
\includegraphics[scale=0.5]{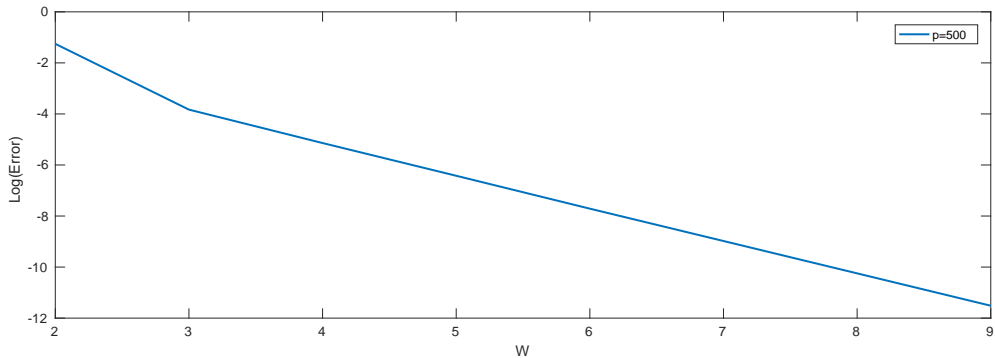}

\caption{Log of relative error against $W$}

\end{figure}

\subsubsection*{Conclusions }

The proposed spectral element method for elliptic interface problem
with nonsmooth solutions is nonconforming and exponentially accurate.
The interface conditions are incorporated as jumps across the interfaces
in appropriate Sobolev norms in the least-squares formulation. The
numerical method is also applicable for general polygonal domains.
The numerical solution has been obtained efficiently and inexpensively
using PCGM. A decoupled block diagonal preconditioner has been used.
More efficient preconditioner for the interface problems is under
investigation.

\end{document}